\documentclass[12pt,a4paper]{amsart}
\usepackage{amssymb,amsmath,amsthm}
\usepackage{graphicx}
\usepackage{epsfig}
\usepackage{tikz}

\long\def\onefigure#1#2{%  #1 picture,  #2  caption
\begin{figure*}[tbp]
\begin{center}
#1
\end{center}
\caption{#2}
\end{figure*}
} %end onefigure def

\newcommand{\lipefig}[2]  % labeled Ipe figure
{\onefigure{\mbox{\psfig{file=#1.eps}}}{\label{f:#1} #2} }

\newtheorem{theorem}{Theorem}[section]
\newtheorem{lemma}{Lemma}[section]
\newtheorem{corollary}{Corollary}[section]

\newtheorem{claim}{Claim}[section]

\newcommand{\de}{\delta}

\newcommand{\eps}{\varepsilon}

\newcommand{\conv}{\mathrm{conv}}

\newcommand{\cone}{\mathrm{cone}}
\newcommand{\prob}{\mathrm{Prob}}

\numberwithin{equation}{section}

\begin{document}

\title{The cocked hat}
\author{Imre B\'ar\'any, William Steiger, Sivan Toledo}

\begin{abstract} We revisit the cocked hat -- an old problem from navigation -- and examine under what conditions its old solution is valid.
\end{abstract}

\maketitle

\section{Introduction}\label{sec:introd}

\begin{figure}
\centering
\includegraphics[scale=0.8]{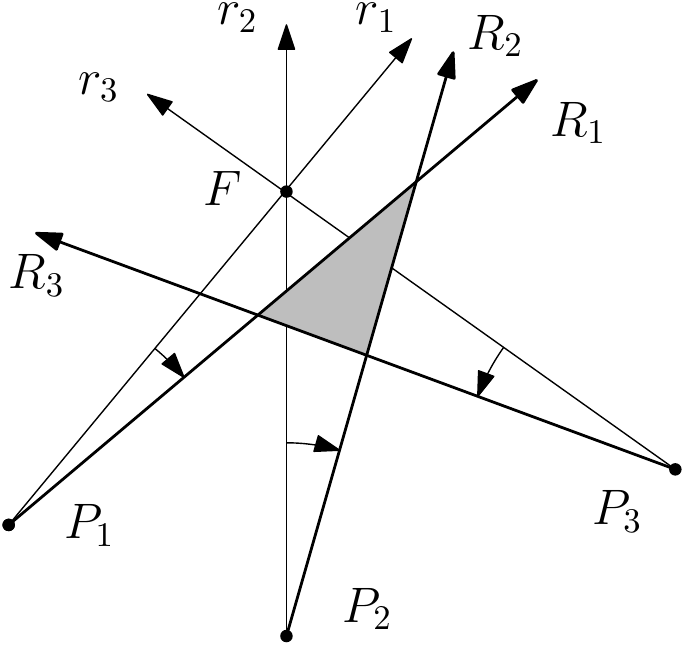}
\caption{Three observation points $P_1,P_2,P_3$, the target $F$, the three rays, and the cocked hat (shaded).}
\label{fig:CockH}
\end{figure}

Navigators used to plot on a map \emph{lines of position} or \emph{lines
of bearing}, which are rays emanating from a landmark (e.g., a lighthouse
or radio beacon) at a particular bearing (angle relative to north)
that was estimated to be the direction from the landmark
(which we also refer to as observation point) to the ship
or plane. Two such rays usually intersect at a point, which the navigator
would take as an estimate of the true position of the craft. Navigators
were encouraged to plot three rays, to make position estimation more
robust. The three rays normally created a triangle, called a \emph{cocked
hat}~\cite{OxfordCompanionToShipSea}, as shown in Figure~\ref{fig:CockH}.
The properties of the cocked hat were investigated thoroughly~\cite{anderson_1952,cook_williams_1993,cotter_1961,
daniels_1951,AdmiraltyNavigationManual3_1938,stansfield_1947,stuart_2019,williams_1991},
to help navigators interpret it and make good navigation decisions.
The aim of this paper is to analyze the conditions under which an
elegant property of the cocked hat holds. That property had been stated
without a proof more than 80 years ago~\cite{AdmiraltyNavigationManual3_1938},
proved informally (and essentially incorrectly) 70 years ago~\cite{stansfield_1947},
and has been widely disseminated ever since~\cite{cook_williams_1993,daniels_1951,TheScienceOfNavigation,
OUVideo_cocked_hat,williams_1991},
including in course material~\cite{OUVideo_cocked_hat} and in a
popular science book~\cite{TheScienceOfNavigation}.

The property that we are interested in is the probability of the cocked
hat containing the true position being $1/4$. Under what conditions
is this statement true?

This claim first appeared in a 1938 navigation manual~\cite[page~166]{AdmiraltyNavigationManual3_1938},
without a proof and with only informal conditions on the error angles
at the three landmarks, which we denote $P_{1}$, $P_{2}$, and $P_{3}$ (see Figure~\ref{fig:CockH}).
The error angles $\epsilon_{1}$, $\epsilon_{2}$, and $\epsilon_{3}$
are between the plotted rays, which we denote $R_{1}$, $R_{2}$,
and $R_{3}$ and the rays $r_i$ from $P_{i}$ to the true position of the
craft, which we denote by $F$ (see Figure~\ref{fig:CockH}).
The informal conditions are that the errors
are independent (the manual does not use this term, but this is what
it means) and fairly small, around 1 degree. A 1947 article by Stansfield~\cite{stansfield_1947}\footnote{Stansfield developed
the results published in the paper while serving
in Operational Research Sections attached to the Royal Air Force Fighter
Command and Coastal Command during World War II.} cites the claim,
gives more formal conditions for it, and sketches
a proof. The conditions that Stansfield specified are remarkably weak:
he claims that the result would hold if only two of the three errors have zero
median. Stansfield writes that this assumption is equivalent to the following:
``for two of the stations the observed bearings are equally likely to pass to
the right or the left of the true position''. A 1951 article by Daniels~\cite{daniels_1951}
states Stansfield's result in a more modern statistical language,
saying that the cocked hat is a $25\%$ distribution-free confidence
region; the term \emph{distribution free} means that the result is
not dependent on a particular error distribution, say Gaussian, but
only on a parameter of the distribution, here the zero median\footnote{Daniels
was a statistician and served as the president of the Royal
Statistical Society from 1974 to 1975. His paper incorrectly states
that the Admiralty Navigation Manual proves the $25\%$-probability
result; it does not; the first proof sketch appears in Stansfield's
paper. }. Daniels then considers the case of $n$ landmarks and $n$ rays
starting from there. The lines of these rays split the plane into finitely connected
components, some of them bounded, some of them not. Daniels claims without
proof a particular formula, $\frac {2n}{2^n}$, for the probability that $F$ belongs
to the union of the unbounded components. The $25\%$-probability result
was incorrectly extended again by Williams\footnote{Williams was
a professional air navigator and served as president
of the Royal Institute of Navigation from 1984-1987~\cite{charnley_1993}.}
in 1991. He claimed specific probabilities that the open regions
around the cocked hat contain $F$, again with only an informal specification
of the assumptions and with only a sketch of the proof. Williams's
claims were shown to be false by Cook~\cite{cook_williams_1993},
using specific error distributions to which Williams answered with a witty
(but scientifically wrong) rebuttal. Cook also repeated the claim that the
probability of the cocked hat contains $F$ is $1/4$.

Our aim in this paper is to show that the $25\%$-probability result
is valid only for error distributions that guarantee that the
three rays intersect at three distinct points and form a triangle.

We note that the use of the cocked hat in navigation is today obsolete, having been replaced by estimation of confidence regions, usually circles or ellipses, by computer algorithms.

\section{Generalizations to rays that do not intersect}

Two rays in the plane can intersect, but they can also fail to intersect.
Lines of position plotted by navigators almost always intersected,
because the error angles were small. Also, navigators were taught
to choose landmarks so that no angle at the intersection is smaller
than about 50 degrees -- a small angle at the intersection implies
ill conditioning (high sensitivity of the intersection point to bearing
errors).

Stansfield's formulation of the problem uses much more general assumptions
on the errors, and no assumption about angles at the intersections.
Stansfield, Daniels, and the authors that followed only require that
the three errors $\eps_1,\eps_2,\eps_3 \in(-\pi,\pi]$ are random, independent,
and that the median of their distributions is zero. We replace the
zero-median assumption by a consistent but slightly more general condition, namely

\begin{equation}\label{eq:median}
\prob(\epsilon_{i}<0)=\prob(\epsilon_{i}>0)=\frac 12
\end{equation}
\noindent for every $i \in [n]$, where $[n]$ is a shorthand for the set
$\{1,2,\ldots,n\}$.
This means in particular that the target is {\sl never on} $R_i$
which is a necessity because if
$\prob(\epsilon_{i}=0)>0$ were allowed, then $\prob(F \in \Delta)$
could be close to one (e.g., if $\prob(\epsilon_{i}=0)$
is close to one), implying that the $1/4$ result does not hold in this case.
We also consider the restriction of the errors to $[-\pi/2,\pi/2]$.

Under these weak assumptions on the error distribution, the three
rays might fail to form a triangle (the cocked hat). How can we formally
express the $25\%$-probability result when rays may fail to intersect?
We propose four ways to express the result; the first three are fairly
natural but are not sufficient for the $1/4$ result, even under the restriction
$\epsilon_{i}\in[-\pi/2,\pi/2]$; the fourth is not particularly natural
but is the only correct statement of the result.

\smallskip
{\bf Conjunction formulation.} The probability that the three
rays intersect at three points and that the triangle that they form
contains $F$ is $1/4$. In this formulation, we allow error distributions
that could generate non-intersecting rays and we hope to prove that
the probability that the rays intersect at fewer than three points
or that the triangle does not contain $F$ is exactly $3/4$. This is false.

\smallskip
{\bf Conditional probability formulation.} The conditional probability
that the triangle that the rays form contains $F$, conditioned on
the rays forming a triangle, is $1/4$. In this formulation we again
allow error distributions that generate non-intersecting rays,
and we hope to prove that if the rays intersect at three points, then
the probability that the triangle contains $F$ is $1/4$. We do not
care with what probability the rays fail to form a triangle.
This again is false.

\smallskip
{\bf Lines formulation.} We extend the rays $r_{i}$ to infinite
lines $\ell_{i}$, which always form a triangle, and we hope to show
that the triangle that they form contains $F$
with probability $1/4$. Here we must restrict $\epsilon_{i}\in[-\pi/2,\pi/2]$,
otherwise the same line could appear both on the left and on the right
of $F$. Again, this claim is false.

\smallskip
{\bf Constrained distribution formulation.} We assume that the
distribution of errors is such that every pair of rays always intersects
and we hope to show that the probability that the
triangle contains $F$ is $1/4$. We do not permit distributions under
which two of the rays might fail to intersect. We show below that in
this case $\prob(F \in \Delta)=\frac 14$.

\smallskip
We note that from the navigator's perspective, the conditional probability
is the most natural. You plot three rays. If they do not intersect
at three points, you discard the measurements and try again, because
you either picked bad observation points (e.g., two of them and your ship are almost
collinear) or at least one of the bearings is way off. If they do
intersect at three points, you want to know the (conditional) probability
that the cocked hat contains $F$. From the statistician's perspective,
any of the first three formulations makes sense. The fourth makes
less statistical sense, because it is unusual to assume that independent
error distributions satisfy some global structural constraint. In
particular, it appears that Daniels may have believed that the lines
formulation is correct, because he writes about geometrical lines
in the plane, not about rays. He writes ``a particular set of $n$
lines, no two of which are parallel, divides the plane in to $\frac{1}{2}(n^{2}+n+2)$
polygons''.

\section{Counterexamples}\label{sec:counterex}

We now show that the Conjunction, Conditional probability, and Lines
formulation are all false by giving counterexamples. Every example is a {\sl two-ray distribution}
that is concentrated on two rays $R_i^+$ and $R_i^-$:
$\prob(R_i=R_i^+)=\prob(R_i=R_i^-)=\frac 12$. This is no coincidence as we will see
at the end of this section.

\begin{figure}
\centering
\includegraphics[scale=0.7]{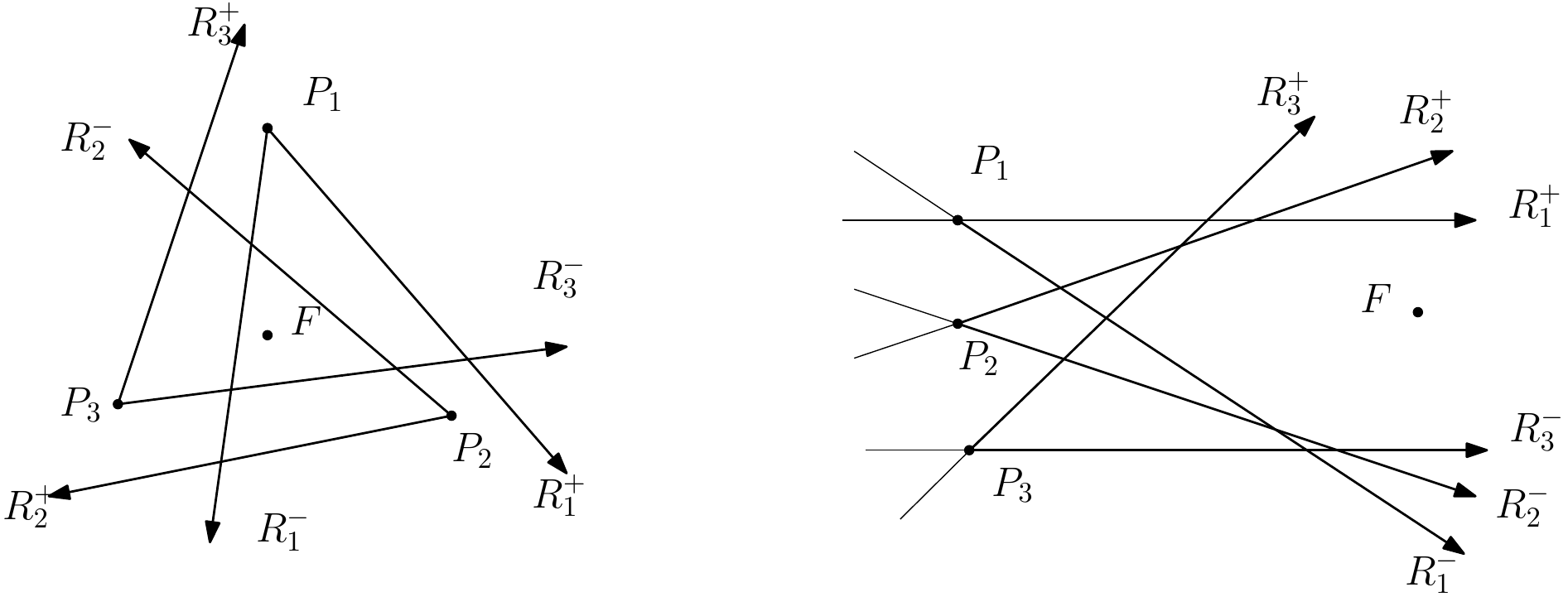}
\caption{Two counterexamples.}
\label{fig:CExamp1}
\end{figure}

In the first example $F$ is in the centroid of an equilateral triangle whose vertices
are $P_{1}$, $P_{2}$, and $P_{3}$. Figure~\ref{fig:CExamp1}
(left) shows the two-ray error distributions. It is easy to see that $R_{i}^{+}$
intersects neither $R_{i+1}^{+}$ nor $R_{i+1}^{-}$ (subscripts are meant modulo $3$).
Therefore, if $R_{i}^{+}$ is selected, then a
cocked hat does not form. On the other hand, if $R_{1}^{-}$, $R_{2}^{-}$,
and $R_{3}^{-}$ are selected, then they form a cocked hat that contains
$F$. Therefore,
\begin{eqnarray*}
\prob\left(F\in\triangle|\text{ the rays form a cocked hat }\triangle\right) & = & 1\\
\prob\left(\text{the rays form a cocked hat }\triangle\text{ and }F\in\triangle\right) & = & \frac{1}{8}\;.
\end{eqnarray*}
This shows that both the Conjunction formulation is false and that
the Conditional probability formulation is false. Note that all the
error angles have magnitude less than $\pi/2$, so these formulations
are false even with this restriction.

Figure~\ref{fig:CExamp1} (right) shows another two-ray distribution. The error magnitudes are less than
$\pi/2$, actually as small as you wish.
The true position $F$ lies outside all the triangles that the lines
form, so the probability that the cocked hat (in the Lines formulation)
contains $F$ is zero. We can move $F$ to the right by any amount
and $F \notin \Delta$ will still hold. This example also shows that the conditional
probability that a cocked hat formed by $3$ rays contains $F$ can
also be zero.

\begin{figure}[b!]
\centering
\includegraphics[scale=0.8]{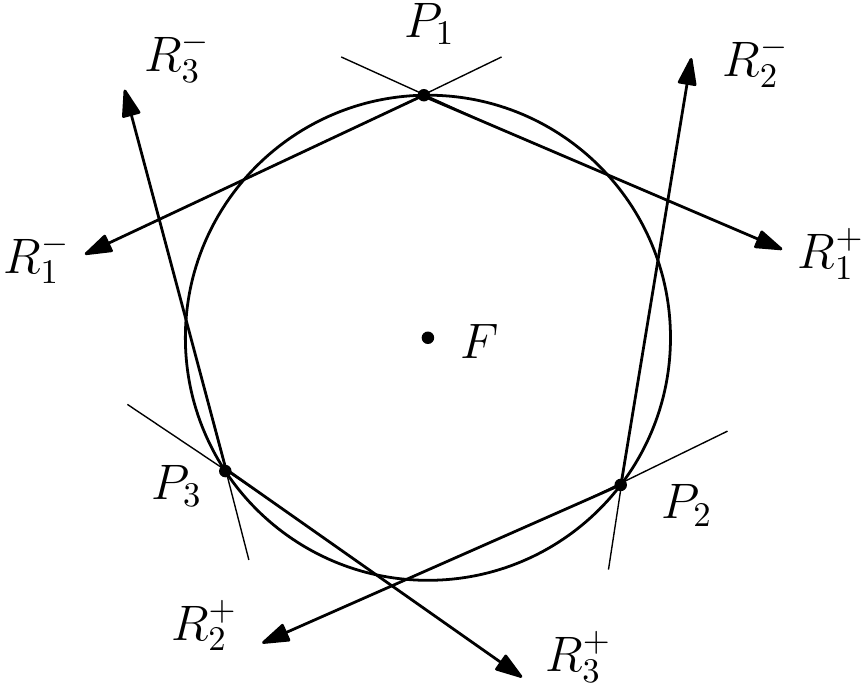}
\caption{The third counterexample.}
\label{fig:CExamp2}
\end{figure}

The last example, given in Figure~\ref{fig:CExamp2}, shows
that the probability that the triangle formed by the extension of
the rays to lines contains $F$ can be $1$. We again note that the
error angles are bounded in magnitude by $\pi/2$.
In this example the three rays do not have three intersection points,
so the cocked hat appears with probability zero. So this is another
counterexample to the Conjunction formulation.

We close this section with a remark on two-ray distributions.
The set of (Borel) probability distributions satisfying
condition (\ref{eq:median}) is convex, and its extreme points are exactly the
two-ray distributions, as one can easily check. Moreover $\prob(F \in \Delta)$
is a linear function on the product of the distributions $\mu_1,\mu_2,\mu_3$
where $\mu_i$ is the probability distribution of the ray $R_i$. Indeed, denoting by $I(E)$ the indicator function of an event $E$, we have
\begin{equation}\label{eq:tworay}
\prob(F \in \Delta)=\int I(F \in \Delta)d\mu_1d\mu_2d\mu_3,
\end{equation}
a linear function of each $\mu_i$, so if it takes the value $\frac 14$
on the two-ray distributions, then it takes the same value on all distributions
satisfying (\ref{eq:median}). We will come back to such distributions in Section~\ref{sec:daniels} again.

\section{Intersecting rays}

We now start the analysis when rays must intersect in pairs. We assume
throughout that the $n+1$ points $P_{1},\ldots,P_{n},F$ are
in general position, so that no three are collinear and so that no
other degeneracies arise.

We introduce some notation. We let $\overrightarrow{XY}$
denote the ray emanating from $X$ in the direction of
$Y$ when $X,Y$ are distinct points in the plane; here we assume that $X \notin \overrightarrow{XY}$.
Thus $r_i=\overrightarrow{P_{i}F}$ is the ray starting at
$P_{i}$ in the direction of the target $F$, and $\ell_i$ is the line containing $r_i$.
From each $P_i$ out goes a random ray $R_i$ making a (signed) angle $\eps_i\in (-\pi,\pi)$
with $r_i$. Our basic assumption, besides (\ref{eq:median}), is
 that two random rays always intersect that is for distinct $i,j \in  [n]$
\begin{equation}\label{eq:nempty}
\prob(R_i\cap R_j=\emptyset)=0.
\end{equation}
So ray $R_i$ and $R_j$ intersect almost surely but their intersection point
is not $P_i$ or $P_j$ because of  our convention
that $X \notin \overrightarrow{XY}$ .

Further notations: $h_i^-$ resp. $h_i^+$ are the halfplanes bounded by $\ell_i$ with $h_i^-$
consisting of points $X$ such that the ray $\overrightarrow{P_iX}$ comes from a clockwise rotation from $r_i$ with angle less than $\pi$, and $h_i^+$ is its complementary halfplane.
When $r,r'$ are two rays we denote
by $\cone(P_i,r,r')$ the cone whose apex is $P_i$ and whose bounding rays are translated copies
of $r$ and $r'$. Such a cone always has angle less than $\pi$, because $r$ and $r'$ will never have
opposite directions.

Define $C_{ij}=\cone(P_i,r_j,\overrightarrow{P_iP_j})$ for distinct $i,j \in [n]$.

\begin{lemma}\label{l:cones} The cone $C_{ij}$ contains $r_i$ and $\prob( R_i \subset C_{ij})=1$.
\end{lemma}

\begin{figure}
\centering
\includegraphics[scale=0.75]{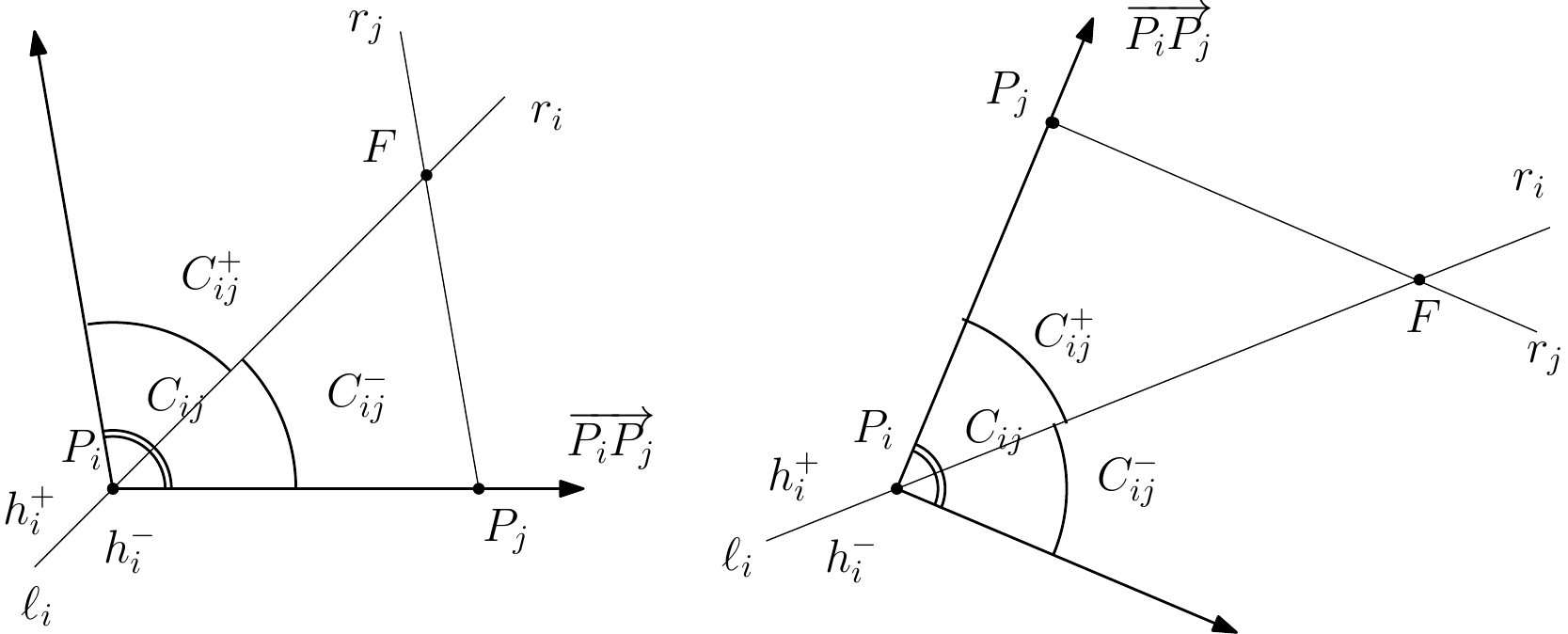}
\caption{Illustration for Lemma~\ref{l:cones}: the case $P_j \in h_i^-$ on the left, the case $P_j \in h_i^+$ on the right.}
\label{fig:cones}
\end{figure}

{\bf Proof.} Assume first that $P_j \in h_i^-$. We define first the cones $C_{ij}^-=\cone(P_i,r_i, \overrightarrow{P_iP_j})$ and $C_{ij}^+=\cone(P_i,r_i,r_j)$, see Figure~\ref{fig:cones}. Note that the angle of $C_{ij}^-$ (resp. $C_{ij}^+$) is equal to the angle at $P_i$ (and at $F$) of the triangle with vertices $P_i,P_j,F$. Then $C_{ij}=C_{ij}^- \cup C_{ij}^+$ because the angle of this cone is the sum of the angles of $C_{ij}^-$ and $C_{ij}^+$ so smaller than $\pi$. Then $r_i \subset C_{ij}$ indeed as shown in  Figure~\ref{fig:cones}, left.

Suppose now that $\eps_i>0$ which is the same as $R_i \subset h_i^+$. If $R_i$ does not lie in $C_{ij}^+$, then $R_i \subset h_i^+\setminus C_{ij}^+$. The last set is a convex cone, disjoint from $h_j^-$, as they are separated by the line $\ell_j$. So no $R_j$ with $\eps_j<0$ can intersect $R_i$ contradicting (\ref{eq:nempty}). So $R_i \subset C_{ij}^+$.

Let $h$ denote the halfplane containing $F$ and bounded by the line through $P_i$ and $P_j$. Observe that by the previous argument $R_j \subset h$ because the complementary halfplane to $h$ is disjoint from $C_{ij}^+$, so $R_j$ can intersect $R_i \subset C_{ij}^+$ only if it lies in $h$.

Suppose next that $\eps_i<0$. We show that $R_i \subset C_{ij}^-$. If not, then $R_i \subset h_i^-\setminus C_{ij}^-$. The last set is a convex cone again, disjoint from $h$, so $R_i\cap R_j=\emptyset$ for all $R_j$ with $\eps_j<0$ contradicting (\ref{eq:nempty}).

The argument for the case $P_j \in h_i^+$ is symmetric (see Figure~\ref{fig:cones} right) but otherwise identical and is therefore omitted.\qed

\medskip
We remark here that Lemma~\ref{l:cones} implies that the cone $\bigcap_{j\ne i} C_{ij}$ is convex (that is, its angle is smaller than $\pi$), it contains $r_i$, and $\prob(R_i \subset \bigcap_{j\ne i} C_{ij})=1$, of course only if $n\ge 2$. (For $n=1$ condition (\ref{eq:nempty}) is void.) Define $K_i$ as the smallest (with respect to inclusion) convex cone satisfying $\prob(R_i \subset K_i)=1$. Note that $K_i \subset \bigcap_{j\ne i} C_{ij}$. For later reference we state the following corollary.

\begin{corollary}\label{cor:cone} Under conditions (\ref{eq:median}) and (\ref{eq:nempty}) $K_i$ is a convex cone, $r_i\subset K_i$ and $\prob(R_i\subset K_i)=1$ for every $i\in [n]$.
\end{corollary}

\begin{theorem}\label{th:1/4} Under conditions (\ref{eq:median}) and (\ref{eq:nempty})
\[
\prob(F \in \Delta)=\frac 14.
\]
\end{theorem}

{\bf Proof.} Set $T=\conv \{P_1,P_2,P_3,F\}$, the convex hull of $P_1,P_2,P_3$, and $F$. We will have to consider three cases separately: when $T$ is a triangle with $F$ inside $T$ (Case 1), when $T$ is a triangle with $F$ a vertex of $T$ (Case 2), and when $T$ is a quadrilateral (Case 3).

\begin{figure}[t!]
\centering
\includegraphics[scale=0.7]{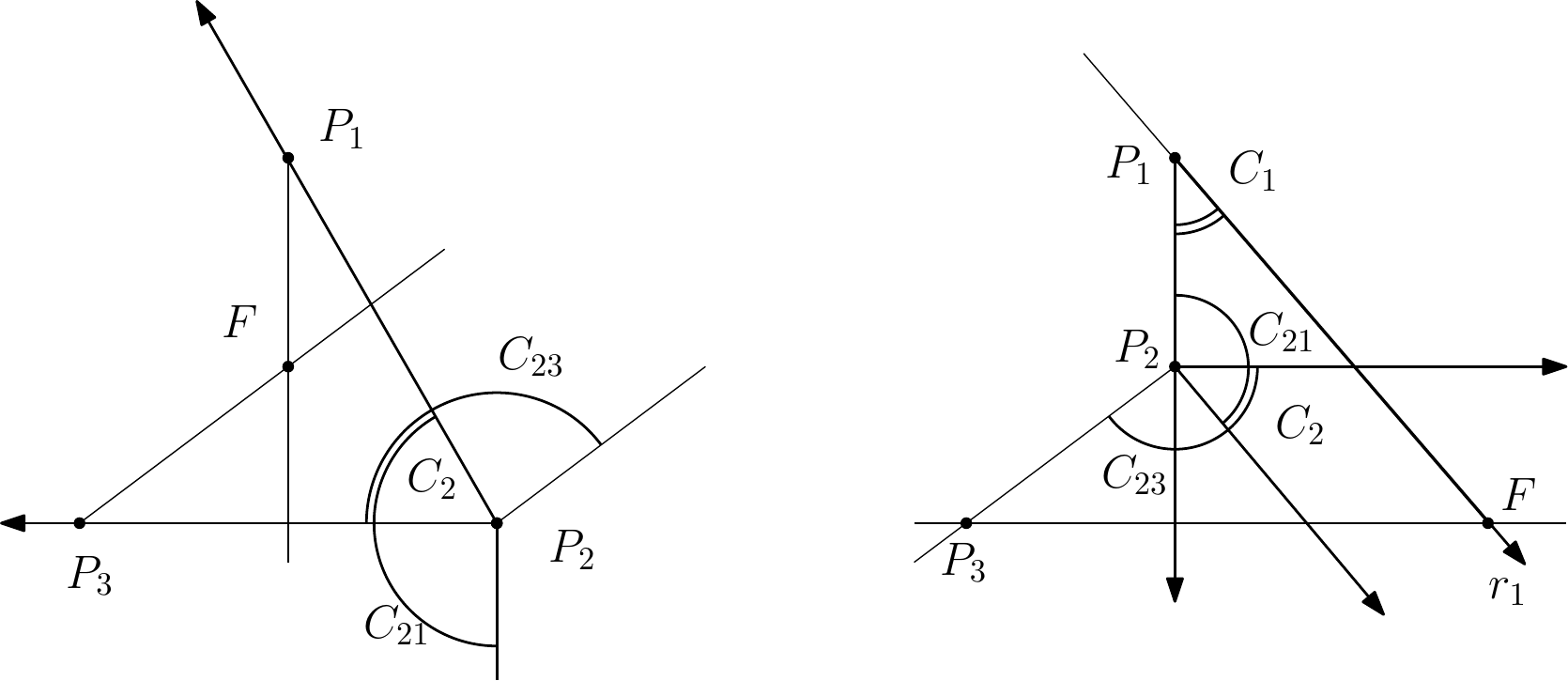}
\caption{The cone $C_2$ in Case 1 (left) and 2 (right).}
\label{fig:onequarter}
\end{figure}

{\bf Case 1.} Define $C_i=\cone(P_i,\overrightarrow{P_iP_{i-1}},\overrightarrow{P_iP_{i+1}})$ for $i=1,2,3$ where the subscripts are taken mod $3$, see Figure~\ref{fig:onequarter} left.

We claim that $R_i \subset C_i$ for all $i$. By symmetry it suffices to show this for $i=2$. By Lemma~\ref{l:cones} $R_2 \subset C_{21}\cap C_{23}$. So it is enough to check that $C_2 = C_{21}\cap C_{23}$, and this is evident: the rays bounding $C_2$ are $\overrightarrow{P_2P_{1}}$ (which bounds $C_{21}$) and $\overrightarrow{P_2P_3}$ (which bounds $C_{23}$).

\begin{figure}[t!]
\centering
\includegraphics[scale=0.8]{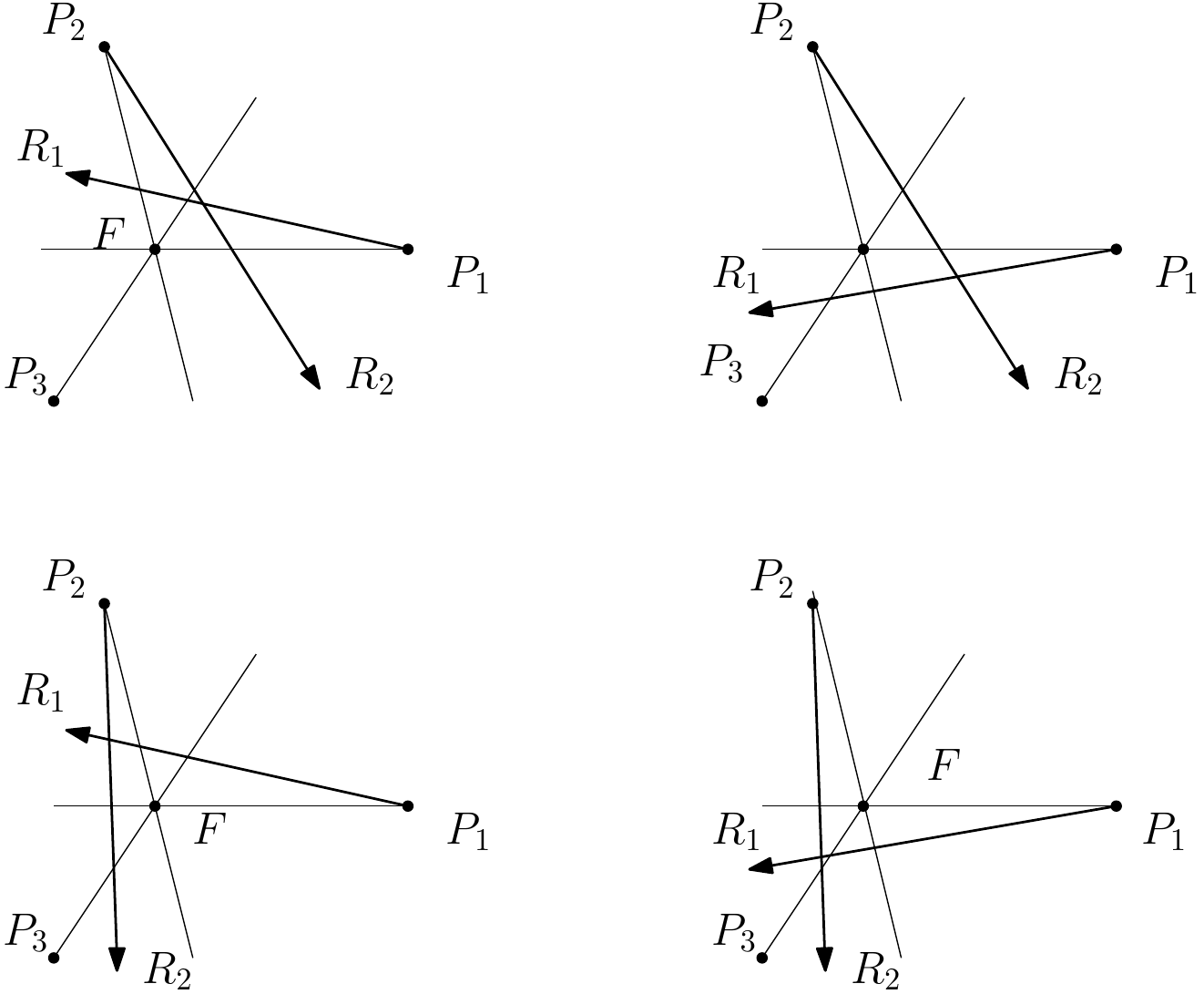}
\caption{Illustration for the proof of Theorem~\ref{th:1/4}.}
\label{fig:8cases}
\end{figure}

We can now finish the proof of the theorem in Case 1. There are $8$
sub-cases with equal probabilities that correspond to the signs of
$\epsilon_{1}$, $\epsilon_{2}$, and $\epsilon_{3}$, as shown in
Figure~\ref{fig:8cases}. Only in two of them, namely when
all $\epsilon_{i}$s have the same sign, we have
$F\in\triangle$, so the probability of this event is $1/4$.

{\bf Case 2.} We assume (by symmetry) that $P_2$ is inside the triangle $T$.
We define the cones $C_1=\cone(P_1,r_2,\overrightarrow{P_1P_2})$, $C_2=\cone(P_2,r_1,r_3)$, and
$C_3=\cone(P_3,r_2,\overrightarrow{P_3P_2})$ and we claim that $R_i \subset C_i$ for all $i$. From Lemma~\ref{l:cones} we have that $R_2 \subset C_{21}\cap C_{23}$. The bounding rays of $C_2$ are a translate of $r_2$ (bounding $C_{21}$) and a translate of $r_3$ (bounding $C_{23}$), so $C_2 = C_{21}\cap C_{23}$ (see Figure~\ref{fig:onequarter} right).

The cases $i=1$ and $3$ are symmetric and very simple. We only consider $i=1$. Again, by Lemma~\ref{l:cones} $R_1 \subset C_{12}$ and then $C_1= C_{12}$ implying $R_1 \subset  C_1$.

Again there are 8 subcases, corresponding to the 8 possible sign patterns of $\eps_1,\eps_2,\eps_3$. It is easy to see that $F \in \Delta$ in exactly two of them.

{\bf Case 3.} We assume again by symmetry that the segment $P_2F$ is a diagonal of the quadrilateral $T$. Define cones $C_1=\cone(P_1,r_2,\overrightarrow{P_1P_3})$, $C_2=\cone(P_2,r_1,r_3)$, and
$C_3=\cone(P_3,r_2,\overrightarrow{P_3P_2})$. We claim again that $R_i \subset C_i$ for all $i$. The proof is similar to the previous ones using Lemma~\ref{l:cones} and is omitted here. Again, $F \in \Delta$ in exactly two  out of the 8 cases.\qed

\section{Daniels' statement}\label{sec:daniels}

We assume now that there are $n\ge 3$ observation points $P_1,\ldots,P_n$ plus the target point $F$ and that these $n+1$ points are in general position. A random ray $R_i$ starts at each $P_i$ satisfying conditions (\ref{eq:median}) and  (\ref{eq:nempty}). The lines of the rays $R_i$ split the plane into connected components, let $U$ denote the union of the $2n$ unbounded components. Here comes Daniels' statement.

\begin{theorem}\label{th:daniels} Under conditions (\ref{eq:median}) and (\ref{eq:nempty})
\[
\prob(F \in U)=\frac {2n}{2^n}.
\]
\end{theorem}

The case $n=2$ is trivial and not interesting. The case $n=3$ is just Theorem~\ref{th:1/4}. We note that condition (\ref{eq:nempty}) is a necessity, even for $n=3$ as the counterexamples in Section~\ref{sec:counterex} show.

We are going to prove this theorem under the assumption that each $R_i$ is a two-ray distribution, that is,
$\prob(R_i=R_i^+)=\prob(R_i=R_i^-)=\frac 12$ and explain, after the proof, how this special case implies the theorem. We also assume that the $2n$ rays $R_i^+,R_i^-$, together with the points $P_1,\ldots,P_n,F$ are in general position. This is not a serious restriction because the general case of two-ray distributions follows from this by a routine limiting argument.

\begin{figure}
\centering
\includegraphics[scale=0.8]{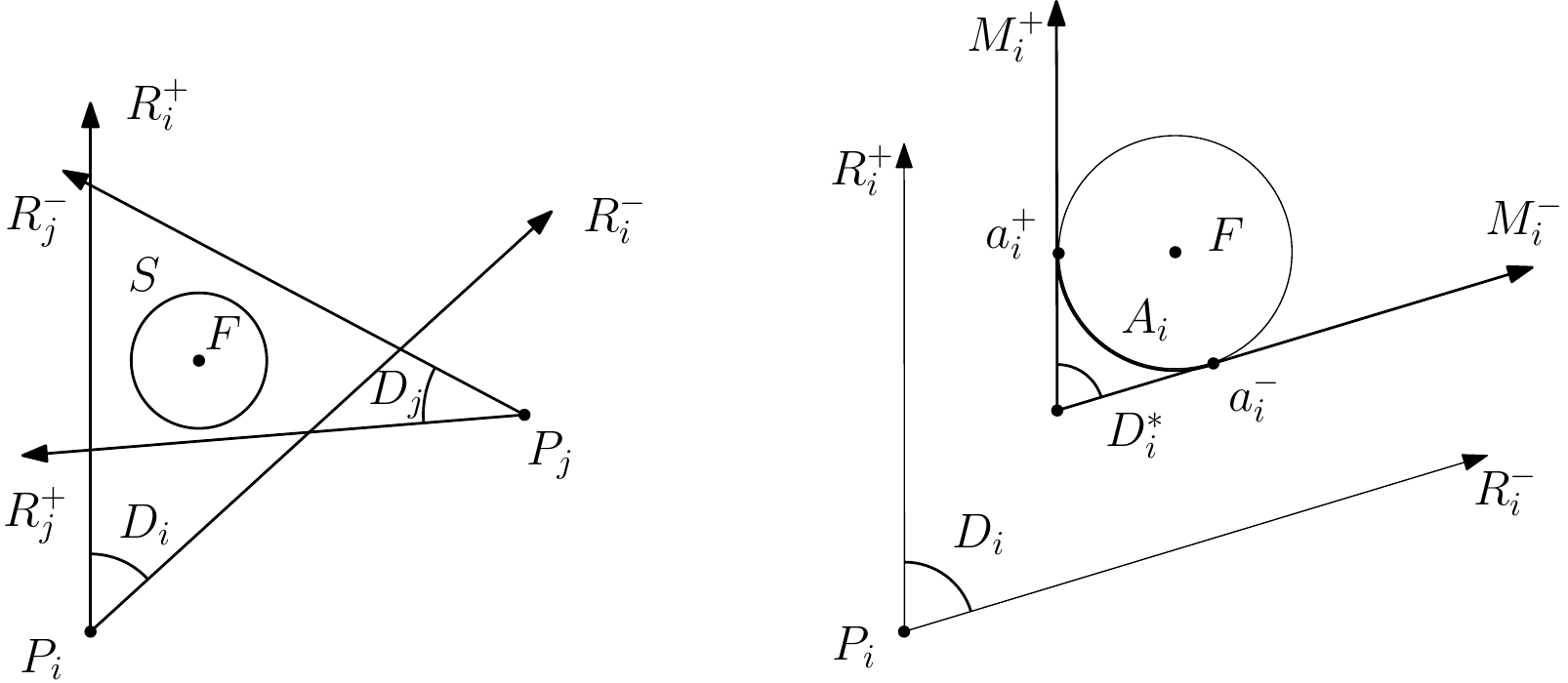}
\caption{The intersection $D_i\cap D_j$ and the translated cone $D_i^*$.}
\label{fig:quadril}
\end{figure}

{\bf Proof.} To simplify the writing, we set $D_i=\cone(P_i,R_i^+,R_i^-)$, which is equivalent to $D_i=\conv(R_i^+ \cup R_i^-)$.  Lemma~\ref{l:cones} implies that $r_i \subset D_i$ for every $i\in [n]$. Let $S$ be a circle centered at $F$ such that $S\subset D_i$ for every $i\in [n]$. Observe that for distinct $i,j\in [n]$, the intersection $D_i\cap D_j$ is a convex quadrilateral containing $S$ and of course $F$, see Figure~\ref{fig:quadril} left. This follows from condition (\ref{eq:nempty}): both $R_i^+$ and $R_i^-$ intersect both  $R_j^+$ and $R_j^-$ and the four intersection points are the vertices of $D_i\cap D_j$ which is then a convex quadrilateral.

Let $L_i^+$ (resp. $L_i^-$) denote the line of the ray $R_i^+$ (and $R_i^-$). For a selection $\de_1,\ldots,\de_n \in \{1,-1\}$ of signs the lines $L_1^{\de_1},\ldots,L_n^{\de_n}$ split the plane into finitely many connected components. We are going to show that out of the $2^n$ possible selections there are exactly $2n$ for which $F$ lies in an unbounded component.

We reduce this statement to another one about arcs on the unit circle. First comes a simpler reduction. Translate each cone $D_i$ into a new (and actually unique) position $D_i^*$ so that its rays touch the circle $S$ (see Figure~\ref{fig:quadril} right). Let $Q_i^+,M_i^+$ (resp. $Q_i^-,M_i^-$) be the translated copies of $R_i^+,L_i^+$ (and $R_i^-,L_i^-$). Note that $D_i^*\cap D_j^*$ is again a convex quadrilateral.

We {\bf claim} next that for a fixed selection $\de_1,\ldots,\de_n$ of signs,  $F$ lies in an unbounded component for the lines $L_1^{\de_1},\ldots,L_n^{\de_n}$ if and only if it lies in the corresponding unbounded component for the lines $M_1^{\de_1},\ldots,M_n^{\de_n}$. This is simple. The point $F$ lies in an unbounded component for the lines $L_i^{\de_i}$ if and only if there is a halfline $R$ starting at $F$ and  disjoint from each $L_i ^{\de_i}$ which happens if and only if $R$ is disjoint from the lines $M_i ^{\de_i}$ as well.

Assume now that $S$ is the unit circle. Let $a_i^+$ (resp. $a_i^-$) be the points where $M_i^+$ (and $M_i^-$) touch $S$, and let $A_i$ be the shorter arc on $S$ between $a_i^+$ and $a_i^-$, see Figure ~\ref{fig:quadril} right. It is clear that $a_i^+$ and $a_i^-$ are not opposite points on $S$ so $A_i$ is welldefined. These arcs completely determine $D_i^*$. They satisfy the conditions
\begin{itemize}
\item[(i)] each $A_i$ is shorter than $\pi$, and
\item[(ii)] $A_i\cup A_j$ is an arc in $S$ longer than $\pi$ for all $i,j\in [n], i\ne j$.
\end{itemize}
The latter condition follows from the fact that $D_i^*\cap  D_j^*$ is a convex quadrilateral.

We call a selection $\de_1,\ldots,\de_n$ {\sl special} if it gives an unbounded component containing $F$.
We {\bf claim} that a selection is special if and only if the points $a_1^{\de_1},\ldots,a_n^{\de_n}$ lie on an arc of $S$ shorter than $\pi$. This is also simple. If there is such an arc, call it $I$ and let $Q$ be the centre point of the complementary arc $S \setminus I$. The ray $\overrightarrow{FQ}$ avoids every line $M_i^{\de_i}$. If there is no such arc, then the connected component containing $F$ (and $S$) is bounded as one can check easily. Therefore it suffices to prove the following lemma.

\begin{lemma}\label{l:special} Under the above conditions there are exactly $2n$ special selections.
\end{lemma}

{\bf Proof.} For a special selection $\de=(\de_1,\ldots,\de_n)$ let $I(\de)$ denote the shortest arc on $S$ containing every $a_i^{\de_i}$, $i\in [n]$. Thus $I(\de)$ is the shorter arc between points $a_i^{\de_i}$ and $a_j^{\de_j}$ for some distinct $i,j\in [n]$, and they are the {\sl endpoints} of $I(\de)$.

\begin{claim} Each $a_i^+$ (and $a_i^-$) is the endpoint of $I(\de)$ for exactly two special selections $\de$.
\end{claim}

It suffices to prove this claim since it implies Lemma~\ref{l:special} and then Theorem~\ref{th:daniels} as well.

\begin{figure}
\centering
\includegraphics[scale=0.9]{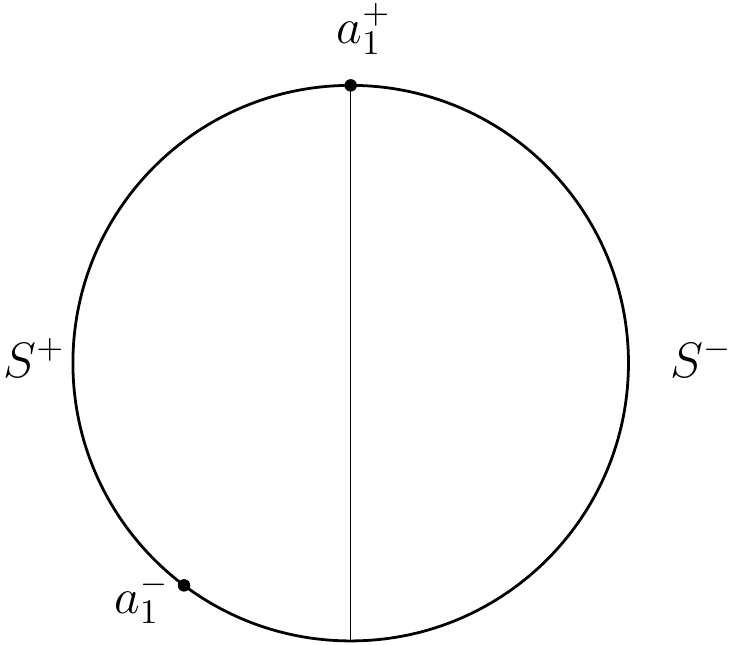}
\caption{The definition of $S^+$ and $S^-$.}
\label{fig:special}
\end{figure}

\smallskip
{\bf Proof} \;of the claim. It is enough to work with $a_1^+$. Using the notation on Figure~\ref{fig:special} we assume that $a_1^-$ is from the halfcircle $S^+$ so $A_1\subset S^+$.

Define $X=\{a_1^+,a_1^-,\ldots,a_n^+,a_n^-\}$ and $Y=X \setminus\{a_1^+,a_1^-\}$. Observe first that $S^+$ can't contain any $A_i$, $i>1$ as otherwise $A_1,A_i\subset S^+$ contradicting (ii). Moreover, $S^-$ can't contain two arcs $A_i,A_j$ with distinct $i,j>1$ because of (ii) again. It follows then that $|S^+\cap Y|=n-1$ or $n-2$.

{\bf Case 1} when $|S^+\cap Y|=n-1$. Then $|S^-\cap Y|=n-1$ as well and $S^+$ contains exactly one element from each pair $\{a_i^+,a_i^-\}$, $i>1$, and then so does $S^-$. This gives exactly two special selection $\de$ and $\eps$ with $I(\de)\subset S^+$ and $I(\eps)\subset S^-$, with $a_1^+$ an endpoint of both.

{\bf Case 2} when $|S^+\cap Y|=n-2$. Then $S^+$ contains no $I(\de)$ with $\de$ special, $|S^-\cap Y|=n$ and so $A_i\subset S^-$ for a unique $i>1$. This gives again two special selections $\de$ and $\eps$ where $a_1^+$ is the endpoint of $I(\de)$ and $I(\eps)$. In fact $\de$ and $\eps$ coincide except at position $i$: $\de_j=\eps_j$ for all $j\in [n]$ but $j=i$ and $\de_1=\eps_1=1$.\qed

\medskip
We explain now how the case of two-ray distributions implies Theorem~\ref{th:daniels}, or rather give a sketch of this and leave the technical details to the interested reader. Assume each ray $R_i$ follows a generic distribution $\mu_i$ for all $i \in [n]$ still satisfying conditions (\ref{eq:median}) and (\ref{eq:nempty}). Note that by Corollary~\ref{cor:cone}, $\prob(R_i \subset K_i)=1$. Using this one can check that every $\mu_i$ can be approximated with high precision by a convex combination of two-ray distributions, each having  $R_i^+,R_i^-\subset K_i$. One has to show as well that this approximation can be chosen so that $R_i^{\de_i}\cap R_j^{\de_j}\ne \emptyset$ for distinct $i,j \in [n]$ and for every choice of signs $\de_i,\de_j$. As in (\ref{eq:tworay}), $\prob(F \in U)$ is a linear function of the underlying distributions $\mu_i$, and this linear function equals $2n/2^n$ on the product of two-ray distributions. Therefore this linear function equals $2n/2^n$ on any convex combination of products of two-ray distributions and consequently $\prob(F \in U)$ must be equal to $2n/2^n$ on the product of the $\mu_i$s.

\bigskip
{\bf Acknowledgements.}  The first author was partially supported by Hungarian National Research grants 131529, 131696, and KKP-133819. The last author was partially supported by grants 965/15, 863/15, and 1919/19 from the Israel Science Foundation and by a grant from the Israeli Ministry of Science and Technology.

\bigskip

\bibliographystyle{plain}
\bibliography{CockedHat}

\bigskip
\noindent
Imre B\'ar\'any \\
R\'enyi Institute of Mathematics,\\
13-15 Re\'altanoda Street, Budapest, 1053 Hungary\\
{\tt barany.imre@renyi.hu}\\
and\\
Department of Mathematics\\
University College London\\
Gower Street, London, WC1E 6BT, UK

\medskip
\noindent
William Steiger\\
Department of Computer Science, Rutgers University\\
110 Frelinghuysen Road,Piscataway, NJ 08854-8019, USA\\
{\tt steiger@cs.rutgers.edu}

\medskip
\noindent
Sivan Toledo\\
Blavatnik School of Computer Science, Tel Aviv University,\\
Tel Aviv 6997801, Israel\\
{\tt stoledo@mail.tau.ac.il}

\end{document}